\documentclass[intlimits]{scrartcl}

\usepackage{mathrsfs}
\usepackage{fancyhdr}
\usepackage[utf8]{inputenc}
\usepackage[T1]{fontenc}
\usepackage{graphics}
\usepackage[pdftex]{graphicx}
\usepackage{amsmath}
\usepackage{amstext}
\usepackage{amsgen}
\usepackage{amsbsy}
\usepackage{amsopn}
\usepackage{amssymb}
\usepackage{amsfonts}
\usepackage{enumerate}
\usepackage{latexsym}
\usepackage{textcomp}

\newcommand{\leertest}[2]{
  \newcommand{\tempa}{}
  \newcommand{\tempb}{#1}
  \ifx\tempa\tempb
     {\ }  
  \else
     \ #2 
  \fi
}
\newcommand{\thmopt}[1]{\ (\upshape{#1})}
\newcommand{\proofopt}[1]{\ {#1}}

\newcounter{theorem}[section]       
\numberwithin{theorem}{section}     
\newcounter{theoremx}[section]      

\newenvironment{Theorem}[1][]
	{
	\refstepcounter{theorem}
	\par\addvspace{0.775\baselineskip}\par
	\mbox{\bfseries\scshape \thetheorem. Theorem\hspace*{-1em}\leertest{#1}{\thmopt{#1}}\hspace*{-0.2334em}:\ }
	\slshape
	}
	{
	\par\addvspace{0.775\baselineskip}\par
	}
\newenvironment{Lemma}[1][]{\refstepcounter{theorem}\par\addvspace{0.775\baselineskip}\par\mbox{\bfseries\scshape \thetheorem. Lemma\hspace*{-1em}\leertest{#1}{\thmopt{#1}}\hspace*{-0.2334em}:\ }\slshape}{\par\addvspace{0.775\baselineskip}\par}
\newenvironment{Cor}[1][]{\refstepcounter{theorem}\par\addvspace{0.775\baselineskip}\par\mbox{\bfseries\scshape \thetheorem. Corollary\hspace*{-1em}\leertest{#1}{\thmopt{#1}}\hspace*{-0.2334em}:\ }\slshape}{\par\addvspace{0.775\baselineskip}\par}
\newenvironment{Prop}[1][]{\refstepcounter{theorem}\par\addvspace{0.775\baselineskip}\par\mbox{\bfseries\scshape \thetheorem. Proposition\hspace*{-1em}\leertest{#1}{\thmopt{#1}}\hspace*{-0.2334em}:\ }\slshape}{\par\addvspace{0.775\baselineskip}\par}

\newenvironment{Theorem-}[1][]{\refstepcounter{theoremx}\par\addvspace{0.775\baselineskip}\par\mbox{\bfseries\scshape \thetheoremx. Theorem\hspace*{-1em}\leertest{#1}{\thmopt{#1}}\hspace*{-0.2334em}:\ }\slshape}{\par\addvspace{0.775\baselineskip}\par}
\newenvironment{Lemma-}[1][]{\refstepcounter{theoremx}\par\addvspace{0.775\baselineskip}\par\mbox{\bfseries\scshape \thetheoremx. Lemma\hspace*{-1em}\leertest{#1}{\thmopt{#1}}\hspace*{-0.2334em}:\ }\slshape}{\par\addvspace{0.775\baselineskip}\par}
\newenvironment{Cor-}[1][]{\refstepcounter{theoremx}\par\addvspace{0.775\baselineskip}\par\mbox{\bfseries\scshape \thetheoremx. Corollary\hspace*{-1em}\leertest{#1}{\thmopt{#1}}\hspace*{-0.2334em}:\ }\slshape}{\par\addvspace{0.775\baselineskip}\par}
\newenvironment{Prop-}[1][]{\refstepcounter{theoremx}\par\addvspace{0.775\baselineskip}\par\mbox{\bfseries\scshape \thetheoremx. Proposition\hspace*{-1em}\leertest{#1}{\thmopt{#1}}\hspace*{-0.2334em}:\ }\slshape}{\par\addvspace{0.775\baselineskip}\par}
\newenvironment{Subl-}[1][]{\refstepcounter{theoremx}\par\addvspace{0.775\baselineskip}\par\mbox{\bfseries\scshape \thetheoremx. Sublemma\hspace*{-1em}\leertest{#1}{\thmopt{#1}}\hspace*{-0.2334em}:\ }\slshape}{\par\addvspace{0.775\baselineskip}\par}
\newenvironment{Theorem*}[1][]{\par\addvspace{0.775\baselineskip}\par\mbox{\bfseries\scshape Theorem\hspace*{-1em}\leertest{#1}{\thmopt{#1}}\hspace*{-0.2334em}:\ }\slshape}{\par\addvspace{0.775\baselineskip}\par}
\newenvironment{Lemma*}[1][]{\par\addvspace{0.775\baselineskip}\par\mbox{\bfseries\scshape Lemma\hspace*{-1em}\leertest{#1}{\thmopt{#1}}\hspace*{-0.2334em}:\ }\slshape}{\par\addvspace{0.775\baselineskip}\par}
\newenvironment{Cor*}[1][]{\par\addvspace{0.775\baselineskip}\par\mbox{\bfseries\scshape Corollary\hspace*{-1em}\leertest{#1}{\thmopt{#1}}\hspace*{-0.2334em}:\ }\slshape}{\par\addvspace{0.775\baselineskip}\par}
\newenvironment{Prop*}[1][]{\par\addvspace{0.775\baselineskip}\par\mbox{\bfseries\scshape Proposition\hspace*{-1em}\leertest{#1}{\thmopt{#1}}\hspace*{-0.2334em}:\ }\slshape}{\par\addvspace{0.775\baselineskip}\par}
\newenvironment{Subl*}[1][]{\par\addvspace{0.775\baselineskip}\par\mbox{\bfseries\scshape Sublemma\hspace*{-1em}\leertest{#1}{\thmopt{#1}}\hspace*{-0.2334em}:\ }\slshape}{\par\addvspace{0.775\baselineskip}\par}

\newenvironment{Def}[1][]{\refstepcounter{theorem}\par\addvspace{0.775\baselineskip}\par\mbox{\bfseries\scshape\thetheorem. Definition\hspace*{-1em}\leertest{#1}{\thmopt{#1}}\hspace*{-0.2334em}:\ }\upshape}{\par\addvspace{0.775\baselineskip}\par}

\newenvironment{Def-}[1][]{\refstepcounter{theoremx}\par\addvspace{0.775\baselineskip}\par\mbox{\bfseries\scshape \thetheoremx. Definition\hspace*{-1em}\leertest{#1}{\thmopt{#1}}\hspace*{-0.2334em}:\ }\upshape}{\par\addvspace{0.775\baselineskip}\par}
\newenvironment{Ex-}[1][]{\refstepcounter{theoremx}\par\addvspace{0.775\baselineskip}\par\mbox{\bfseries\scshape \thetheoremx. Example\hspace*{-1em}\leertest{#1}{\thmopt{#1}}\hspace*{-0.2334em}:\ }\upshape}{\par\addvspace{0.775\baselineskip}\par}
\newenvironment{Rem-}[1][]{\refstepcounter{theoremx}\par\addvspace{0.775\baselineskip}\par\mbox{\bfseries\scshape \thetheoremx. Remark\hspace*{-1em}\leertest{#1}{\thmopt{#1}}\hspace*{-0.2334em}:\ }\upshape}{\par\addvspace{0.775\baselineskip}\par}
\newenvironment{Not-}[1][]{\refstepcounter{theoremx}\par\addvspace{0.775\baselineskip}\par\mbox{\bfseries\scshape \thetheoremx. Notation\hspace*{-1em}\leertest{#1}{\thmopt{#1}}\hspace*{-0.2334em}:\ }\upshape}{\par\addvspace{0.775\baselineskip}\par}
\newenvironment{Rem*}[1][]{\par\addvspace{0.775\baselineskip}\par\mbox{\bfseries\scshape Remark\hspace*{-1em}\leertest{#1}{\thmopt{#1}}\hspace*{-0.2334em}:\ }\upshape}{\par\addvspace{0.775\baselineskip}\par}
\newenvironment{Not*}[1][]{\par\addvspace{0.775\baselineskip}\par\mbox{\bfseries\scshape Notation\hspace*{-1em}\leertest{#1}{\thmopt{#1}}\hspace*{-0.2334em}:\ }\upshape}{\par\addvspace{0.775\baselineskip}\par}

\newenvironment{Proof}[1][]{\par\addvspace{0.775\baselineskip}\par\mbox{\itshape Proof\hspace*{-1em}\leertest{#1}{\proofopt{#1}}\hspace*{-0.2334em}:\ \ }\upshape}{\hfill\mbox{\emph{q. e. d.}}\par\addvspace{0.775\baselineskip}\par}

\newenvironment{Remark}{\par\mbox{\scshape Remark:\ }\upshape}{\par}

\usepackage[austrian,american,british]{babel}
\usepackage{babelbib}
\usepackage{lscape}
\usepackage{tensor}
\usepackage[usenames]{color}

\newcommand{\ee}{_\varepsilon}
\newcommand{\ve}{\varepsilon}

\newcommand{\tens}[1]{\boldsymbol{\mathsf{#1}}}
\newcommand{\vect}[1]{\boldsymbol{#1}}
\newcommand{\upartial}{\partial}

\date{November 16, 2010}
\author{Clemens Hanel \\ \small{University of Vienna, Faculty of Mathematics} \\[-0.8ex] \small{Nordbergstraße 15, 1090 Wien, Austria} \\[-0.8ex] \small{E-Mail: clemens.hanel@univie.ac.at}}
\title{Wave-type equations of low regularity}
\begin{document}

\maketitle

\begin{abstract}
We prove local existence and uniqueness of the Cauchy problem for a large class of tensorial second order linear 
hyperbolic partial differential equations with coefficients of low regularity in a suitable class of
generalized functions.
\end{abstract}

\vspace{1\baselineskip}
\noindent MSC-class: 83C75; 46F30, 35D05, 35Q75

\section{Introduction}\label{Sec:Introduction}

Recently there has been an increasing interest in the wave operator on space times of low regularity, i.e., the metric
d'Alembertian
\[
 \Box_{\tens g}=\sum\limits_{\alpha,\beta=1}^n g^{\alpha\beta}\nabla_\alpha \nabla_\beta=\sqrt{|\det \tens g|}\sum\limits_{\alpha,\beta=1}^n\frac{\partial}{\partial x^\alpha}\left(\sqrt{|\det\tens g|}g^{\alpha\beta}\frac{\partial}{\partial x^\beta}\right)
\]
of a Lorentzian metric $\tens g$ with Levi-Cività connection $\nabla$ of low regularity \cite{Cla98,ViWi00,GMS09}. 
This development draws its physical motivation from an alternative approach to analyse space time singularities 
in general relativity put forward in \cite{Cla96}. There C.~J.~S.~Clarke proposed 
to treat singularities as obstructions rather to the well-posedness of the Cauchy problem for the scalar wave-equation 
than to the extension of geodesics (see e.\,g. \cite[Ch. 8]{HaEl73} for this standard approach).
Physically speaking, a scalar field---which can bee seen as a reasonable replacement for an extended test body, which
is too hard to model in general relativity---is used to detect singularities, leading to the notion of
generalized hyperbolicity.  More precisely, a space time is called generalized hyperbolic and viewed as
``non-singular'' if the scalar wave equation can be uniquely solved locally around each point. Of course,
here one has to invoke a suitable solution concept since the coefficients of the resulting equation 
will generically be of low regularity.

In \cite{Cla98} Clarke proved generalized hyperbolicity of shell crossing singularities using
a suitable weak solution concept. In \cite{ViWi00} Vickers and Wilson proved generalized hyperbolicity
of conical space times using the theory of nonlinear generalized functions of Colombeau \cite{Col84,Col92}.
More precisely, they embedded the conical space time metric (component wise) into the (full) Colombeau algebra 
and proved unique solvability of the wave equation in this framework invoking a refined version of higher order
energy estimates (cf. \cite[Ch. 7]{HaEl73}). They also succeeded in showing that their generalized solution 
has a distributional limit that fits the expectations from physics. Later Grant, Mayerhofer, and
Steinbauer in \cite{GMS09} generalized the work of \cite{ViWi00} to a fairly large class of ``weakly singular'' space times, where essentially the metric was assumed to be locally bounded: Modelling such space time metrics in Colombeau 
generalized functions from the start, they proved generalized hyperbolicity, however, did not relate their
result to more classical notions.

In this paper we consider a related but substantially different problem: We are concerned with the Cauchy problem for
second order linear hyperbolic tensor equations $L\tens u=\tens F$ with low regularity coefficients on a {\em classical} manifold $M$.
It then follows that in any local coordinate system $x^\alpha$ the differential operator $L$ takes the form
\[
	L=g^{\alpha\beta}\upartial_\alpha\upartial_\beta+\text{lower order terms},
\]
where $g^{\alpha\beta}$ are the contravariant components of a Lorentzian metric. In order to write $L$ in an explicit coordinate free way one then introduces a smooth background metric $\tens{\hat g}$ which enables one to write $L$ in the form
\begin{equation}\label{Eq:PDO}
 	(L\tens u)^I_J=g^{ab}\widehat\nabla_a\widehat\nabla_b u^I_J+B^{aIP}_{JQ}\widehat\nabla_a u^Q_P+C^{IP}_{JQ}u^Q_P,
\end{equation}
with $I,J,P,Q$ multiindices (for details see below).

Here $\widehat\nabla$ denotes the Levi-Cività connection with respect to the \emph{smooth} metric $\tens{\hat g}$, and we also work in the framework of non-linear distributional geometry \cite{KuSt02a,KuSt02b}, that 
is generalized functions in the sense of Colombeau \cite{Col84,Col92}.

Our motivation to study this problem is twofold. First, when inspecting the methods used in \cite{ViWi00,GMS09}, 
we find that the (one and only, singular) metric has to play different roles in different places:
as the principal part of the operator, defining the Levi-Cività connection, and defining the main part of
the energy tensor which is the essential tool in deriving the key estimates. Here we 
separate these roles by using the two distinct metrics $\tens g$ and $\tens{\hat g}$, where only $\tens g$, which defines 
the principal part of the operator $L$, is of low regularity. The main benefit of doing so is to gain some new insight into the fine structure of the energy estimates and to improve on questions of regularity of generalized solutions. We remark that an elaborate regularity theory within algebras of generalized functions does exist (see e.\,g.
\cite{Hor04,Obe06a,Obe06b}) and connecting to it seems necessary to relate the results of \cite{GMS09} to more classical
function spaces. In particular, our asymptotic conditions on $\tens g$ are quite different from those of \cite{GMS09} and we also see that the results of \cite{GMS09} can be improved, see Remark \ref{Rem:Compare with GMS09} below.

Secondly, this strategy of using separate metrics essentially parallels the strategy used to derive the
reduced Einstein equations (see \cite[Ch. 7]{HaEl73}), and, in fact, our equation is a linearization of the
reduced Einstein equations. So this work can be seen as a necessary first step to eventually treat these 
quasilinear equations in the generalized functions framework.

This paper is organized in the following way: In the rest of this introduction we fix our notation and
collect some prerequisites from nonlinear distributional geometry. In Section 2 we formulate the Cauchy problem for our class of wave-type equations of low regularity, and we state our main theorem of existence and uniqueness in Section 3. Section 4 is dedicated to the heart of the proof, i.\,e., higher order energy estimates, which we finally finish in Section 5.

Throughout this paper we suppose $M$ to be a separable, smooth, orientable Hausdorff manifold of dimension $n$. Furthermore, we will make use of abstract index notation, see \cite{PeRi84}. In particular, a tensor field $\tens u$ of type $(k,l)$ on $M$, i.\,e., $\tens u:\prod_{i=1}^k\Omega^1(M)\times\prod_{i=1}^l\mathfrak X(M)\to\mathcal C^\infty(\mathbb R)$ is denoted by $u^{i_1\cdots i_k}_{j_1\cdots j_l}$ or for short $u^I_J$ by using multiindices $I$ and $J$ of length $|I|=k$ resp. $|J|=l$ (as usual $\mathfrak X(M)$ denotes vector fields and $\Omega^1(M)$ denotes one-forms). Thus for a vector field $\vect\xi$ we write $\xi^a$ and for a one-form $\vect\omega$ we write $\omega_a$. The tensor product is simply denoted by concatenating the two objects in question, i.\,e., $(\vect\omega\otimes\vect\nu)_{ab}=\omega_a\nu_b$. The operation of tensorial contraction is denoted by using twice the same index letter, e.\,g., $\vect\omega(\vect\xi)=\omega_a\xi^a$. For metric tensors $\tens e$, by a slight abuse of the multiindex notation, we write $e_{IJ}$ for $e_{i_1j_1}\cdots e_{i_kj_k}$, whenever $|I|=|J|$. Furthermore, the inverse of some metric $e_{ab}$ will be denoted by $e^{ab}$ with the same convention in the case of multiindices.
To distinguish abstract index notation from calculations in coordinates, we will always use greek indices for tensorial components in a coordinate system. So, e.\,g., the coordinates of $e^{ab}\omega_a$ will read $\sum_{\alpha}e^{\alpha\beta}\omega_\alpha$.

Now we briefly recall the necessary facts from nonlinear distributional pseudo-Riemannian geometry in the sense of J.-F. Colombeau \cite{Col84,Col92}. For more details see \cite[Sec. 3.2]{GKOS01}.

The key idea of Colombeau generalized functions is regularization of distributions by nets of smooth functions depending on a regularization parameter $\ve\in(0,1]$. The basic definition of Colombeau's (special) algebra on $M$ is

\begin{Def}
	We set $\mathcal E(M):=\mathcal C^\infty(M)^{(0,1]}$, denote compact subsets of $M$ by $K$, and denote by $\mathcal P(M)$ the space of linear differential operators on $M$. Then
	\begin{align*}
		\mathcal E_M(M):= & \{(u\ee)\ee\in\mathcal E(M)|\forall K\,\forall P\in \mathcal P(M)\,\exists N\in\mathbb N:\sup_{p\in K}|Pu\ee(p)|=O(\ve^{-N})\}, \\
		\mathcal N(M):= & \{(u\ee)\ee\in\mathcal E(M)|\forall K\,\forall P\in \mathcal P(M)\,\forall m\in\mathbb N:\sup_{p\in K}|Pu\ee(p)|=O(\ve^{m})\}. \\
	\end{align*}
	The quotient $\mathcal G(M):=\mathcal E_M(M)/\mathcal N(M)$ is the special Colombeau algebra on $M$. Elements in $\mathcal G(M)$ are denoted by $u=[(u\ee)\ee]=(u\ee)\ee+\mathcal N(M)$. Note that generalized functions can be localized in the obvious way, see \cite[Sec. 3.2]{GKOS01}. We introduce generalized numbers as the ring of constants in the special algebra, i.\,e., generalized functions with vanishing derivative. For a characterization using asymptotic estimates we refer to \cite[Sec. 1.2]{GKOS01}.
We obtain generalized tensor fields of type $(k,l)$ by setting
\[
	\mathcal G^k_l(M):=\mathcal G(M)\otimes\mathcal T^k_l(M).
\]
\end{Def}
This allows us to define the notion of a generalized metric on $M$.
\begin{Def}
	A generalized pseudo-Riemannian metric is a symmetric tensor field $\tens g\in\mathcal G^0_2(M)$ such that $\det\tens g$ is invertible in the generalized sense, i.\,e., for any representative $(\det\tens g\ee)\ee$ of $\det\tens g$ we have
	\[
		\forall K\subset M \text{ compact}\,\exists m\in\mathbb N:\inf_{p\in K}|\det(\tens g\ee)|\geq\ve^m.
	\]
\end{Def}
By Theorem 3.2.74 in \cite{GKOS01} on any relatively compact set a generalized pseudo-Riemannian metric possesses a representative $(g\ee)\ee$ consisting of smooth pseudo-Riemannian metrics. This in turn can be used to define the index (see \cite{GKOS01}, Definition 3.2.75) of a generalized metric and finally, we call $\tens g$ a generalized Lorentzian metric if the index equals 1. This fixes the signature of the metric to be $(1,n-1)$ or equivalently $(-,+,+,\dots)$.
Any generalized Lorentzian (even pseudo-Riemannian) metric induces a $\mathcal G(M)$-linear isomorphism from $\mathcal G^1_0(M)$ to $\mathcal G^0_1(M)$. Moreover, we can adopt the usual classification of vector fields into spacelike, timelike, and null by demanding either $\tens g(\vect\xi,\vect\xi)>0$, $\tens g(\vect\xi,\vect\xi)<0$ or $\tens g(\vect\xi,\vect\xi)=0$, where we have used the notion of \emph{strict positivity}: A generalized function $f$ is called strictly positive, denoted by $f>0$, if
\[
	\forall K\subset M \text{ compact}\,\exists m\in\mathbb N:\inf_{p\in K}f\ee\geq\ve^m.
\]
For further details and, in particular, for a pointwise description, we refer to \cite{KuSt02b,May08}.

\section{A low regularity Cauchy problem}\label{Sec:Setting}

For the rest of the paper we fix a smooth Lorentzian manifold $(M,\tens{\hat g})$.
We are interested in the local forward-in-time Cauchy problem for hyperbolic linear partial differential operators with generalized coefficients, i.\,e.,
\begin{equation}\label{Eq:Wave equation}
	(L\tens u)^I_J=g^{ab}\widehat\nabla_a\widehat\nabla_b u^I_J+B^{aIP}_{JQ}\widehat\nabla_a u^Q_P+C^{IP}_{JQ}u^Q_P=F^I_J
\end{equation}
with initial data
\begin{align}\label{Eq:Initial data}
	\tens u|_{\Sigma_0}= & \tens u_0 & \widehat\nabla_{\vect\xi}\tens u|_{\Sigma_0}= & \tens u_1.
\end{align}
Here $\widehat\nabla$ is the Levi-Cività connection of the smooth metric $\tens{\hat g}$ and $\Sigma_0$ denotes some initial surface with normal vector field $\vect\xi$ (to be detailed below). We denote by $\tens g=g^{ab}$, $\tens B=B^{aIP}_{JQ}$ and $\tens C=C^{IP}_{JQ}$ the low regularity coefficients of $L$ to be modelled in $\mathcal G$. In particular, $\tens g$ will be a generalized Lorentzian metric, and $\tens B$, $\tens C$ will be generalized tensor fields of suitable type, subject to additional conditions to be specified later. Also, the data $\tens F$, $\tens u_0$, and $\tens u_1$ are allowed to be generalized, that is $\tens F\in\mathcal G^k_l$ and $\tens u_0$, $\tens u_1\in\mathcal G^k_l(\Sigma_0)$. We then look for solutions $\tens u\in\mathcal G^k_l$ at least locally.

Our first task is to specify a class of generalized Lorentzian metrics suitable to act as a principal part of $L$. So let $\tens g\in\mathcal G_2^0(M)$ be a generalized Lorentzian metric with representative $(\tens g\ee)\ee$. We want to make sure that locally there exists a suitable foliation of $M$. To this end we chose a relatively compact set $U\subset M$ and ask for the existence of a function $h\in\mathcal C^\infty(U)$ such that $\vect\sigma:=\mathrm dh$ is timelike with respect to $\tens g$. This will be implied by the existence of $M_0>0$ such that for all $\ve$
\begin{equation}\label{Eq:Bounds on g_e}
	\frac {1}{M_0}\leq -\tens g\ee^{-1}(\vect\sigma,\vect\sigma)\leq M_0
\end{equation}
on $U$. Indeed the level surfaces $\Sigma_\tau:=\{q\in U|h(q)=\tau\}$ with $\tau\in[0,\gamma]$ for some $\gamma>0$ are spacelike hypersurfaces with respect to $\tens g\ee$, $\Sigma_0$ beeing the initial surface from \eqref{Eq:Initial data}.

Next we specify a number of asymptotic conditions on $\tens g$, $\tens B$, and $\tens C$. To this end we make use of a smooth Riemannian metric $\tens m$ on $M$ to define the \emph{``pointwise'' norm} of a smooth tensor field $\tens v$, i.\,e.,
\[
|\tens v|^2:=m^{IJ}m_{KL}v^K_I v^L_J.
\]
Note that since we work locally our conditions are in fact independent of the choice of $\tens m$. We now suppose
\begin{enumerate}
	\item[(i)] For every representative $\tens g\ee$, $\tens B\ee$, and $\tens C\ee$ of $\tens g$, $\tens B$, resp. $\tens C$, we demand for all $K$ compact in $U$
	\begin{align*}
		\sup_K|\tens g^{-1}\ee|= & O(1) & \sup_K|\widehat\nabla \tens g^{-1}\ee|= & O(1) \\
		\sup_K|\tens g\ee|= & O(1) & \sup_K|\tens B\ee|= & O(1) \\
		\sup_K|\tens C\ee|= & O(1)
	\end{align*}
	as $\ve$ tends to zero (Observe that $\widehat\nabla$ denotes the Levi-Cività connection associated with $\tens{\hat g}$ and not with $\tens g$.).
	\item[(ii)] For any representative $(\tens g\ee)\ee$ on $U$, the level set $\Sigma_0$ is a past compact spacelike hypersurface such that $\upartial I\ee^+(\Sigma_0)=\Sigma_0$, where $I\ee^+(\Sigma_0)\subset U$ denotes the future emission. Moreover, there exists a nonempty open set $A\subset M$ and some $\ve_0>0$ such that $A\subset\bigcap_{\ve\leq\ve_0}\overline{I\ee^+(\Sigma_0)}$.
\end{enumerate}
Some remarks on this conditions are in order: Condition (i) gives $\ve$-independent bounds on the coefficients $\tens g$, $\tens B$, and $\tens C$ necessary later on to control the asymptotic behaviour of the energy integrals. Strictly speaking the coefficients $\tens B$, and $\tens C$ depend on the choice of $\tens{\hat g}$; in fact the solvability of the differential equation does not, since we estimate on compact sets. Observe that condition (i) also implies \eqref{Eq:Bounds on g_e} and $\sup_K|\widehat\nabla\tens g|=O(1)$ for all compact $K$. Condition (ii) is necessary to guarantee the existence of classical solutions on the level of representatives on a common domain. Compared with \cite{GMS09}, we demand less regularity on higher order derivatives of the coefficients but more regularity on the first order derivatives of $\tens g$. The latter becomes necessary since the connection $\widehat\nabla$ is not derived from $\tens g$ and therefore $\widehat\nabla\tens g$ is generally non-vanishing.

To simplify the notation, we introduce the following abbreviations:
\begin{itemize}
	\item $\vect\xi\ee:=\tens g\ee^{-1}(\vect\sigma,\cdot)$, the generalized vector field associated to $
	\vect\sigma$. Note that $\vect\xi\ee$ is indeed a generalized vector field, while $\vect\sigma$ is not. Furthermore, note that $\tens g\ee(\vect\xi\ee,\vect\xi\ee)=\tens g\ee^{-1}(\vect\sigma,\vect\sigma)$. This also gives that $\sup_K|\vect\xi\ee|=O(1)$ since $|\vect\xi\ee|^2=g\ee^{ab}\sigma_a g\ee^{cd}\sigma_cm_{bd}$ and $\tens g\ee$ is $O(1)$.
	\item $\sigma:= (-\tens{\hat g}^{-1}(\vect\sigma,\vect\sigma))^{1/2}$, the norm of $\vect\sigma$ measured in terms of the smooth metric $\tens{\hat g}.$
	\item $\vect{\hat\sigma}:=\vect\sigma/{\sigma}$, the unit normal to the hypersurfaces $\Sigma_\tau$ measured with respect to $\tens{\hat g}$.
\end{itemize}

\section{The main theorem}

Following the general strategy for solving
differential equations on a space of Colombeau generalized functions, we have to
pursue the following tasks: We start by writing out the initial value
problem in terms of representatives, i.\,e.,
\begin{align}\label{Eq:Representative PDE}
L\ee\tens u\ee= & \tens F\ee & \tens u\ee|_{\Sigma_0}= & \tens u_{0,\ve}
& \widehat\nabla_{\vect\xi\ee}\tens u\ee|_{\Sigma_0}= & \tens u_{1,\ve}.
\end{align}
Then, using classical theory, we solve separately for each $\ve$,
obtaining a net $(\tens u\ee)\ee$ which is a candidate for a generalized
solution. However, to obtain existence in generalized functions, we first
have to ensure that the $\tens u\ee$ are defined on a common domain
(which is guaranteed by condition (ii) in Section \ref{Sec:Setting}), and
then we have to prove moderateness of the net $(\tens u\ee)\ee$. To
obtain  uniqueness of solutions, we have to show independence of the
class $[(\tens u\ee)\ee]$ of the choice of representatives of the data
$\tens F$, $\tens u_{0}$, and $\tens u_{1}$. Observe that the latter
statement amounts to proving a stability property of the problem. Note that we do not have to proof independence of the representatives of the coefficients of $L$ since $\mathcal G$ is a differential algebra. In
this way we will provide a proof of
\begin{Theorem}\label{Th:Main Theorem}
	Let $(M,\tens{\hat g})$ be a smooth Lorentzian manifold, and let $L$ be a 2nd order partial differential operator of the form \eqref{Eq:PDO} with coefficients $\tens g$, $\tens B$, $\tens C$, where $\tens g\in\mathcal G^0_2(M)$ is a generalized Lorentzian metric, and $\tens g$, $\tens B$, and $\tens C$ are subject to conditions (i) and (ii) above. Let $\Sigma_0$ be a hypersurface spacelike with respect to $\tens g$, locally described by $h^{-1}(\{0\})$ for $h\in\mathcal C^\infty(U)$. Then for any point $p\in\Sigma_0$, there exists an open neighbourhood $V$ of $p$ such that the initial value problem \eqref{Eq:Wave equation}, \eqref{Eq:Initial data} has a unique solution $\tens u\in\mathcal G^k_l(V)$.
\end{Theorem}

\section{Energy estimates}

In substance the proof of Theorem \ref{Th:Main Theorem} relies on higher order energy estimates performed relatively to the foliation of $U$ into the spacelike (with respect to $\tens g$) hypersurfaces $\Sigma_\tau$, $0\leq\tau\leq\gamma$ of Section \ref{Sec:Setting}. We will relate the energy of the solution $\tens u$ on a surface $\Sigma_\tau$ to the energy on the initial surface $\Sigma_0$. To this end we introduce some more notation (see Figure 1).
\begin{figure}\label{Fig:Foliation}
	\begin{center}
		\includegraphics[width=10.7cm]{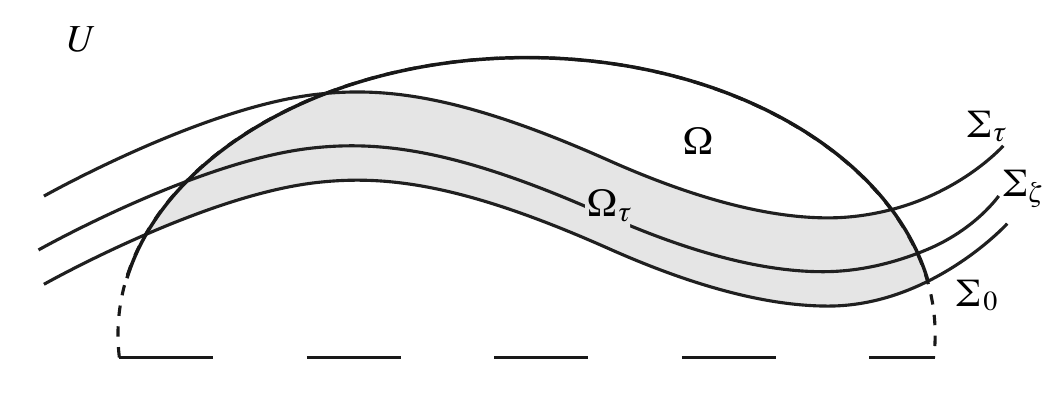}
		\caption{Local foliation of space time}
	\end{center}
\end{figure}
For any set $\Omega\subset U$, we define
\begin{align*}
	\Omega_\tau:= & \Bigl(\Omega\cap\bigcup\nolimits_{0\leq\zeta\leq\tau}\Sigma_\zeta\Bigr)^\circ\text{ and} \\
	S_\tau:= & \Sigma_\tau\cap\Omega. \\
\end{align*}
Now let $p\in\Sigma_0$ and let $\Omega\subset U$ be a relatively compact neighbourhood of $p$ such that $\upartial\Omega$ is spacelike (such a neighbourhood always exists since $\tens g^{-1}$ is locally uniformly bounded).
We denote by $\vect{\hat\mu}$ the volume form on $M$ with respect to $\tens{\hat g}$ and by $\vect{\hat\mu}_\tau$ the induced $(n-1)$-form on $S_\tau$ such that $\vect{\hat\sigma}\wedge\vect{\hat\mu}_\tau=\vect{\hat\mu}$, respectively $\iota_{\vect{\hat\sigma}}\vect{\hat\mu}=\vect{\hat\mu}_\tau$. This allows us to give
\begin{Def}
	Let $\tens v$ be a smooth tensor field, $0\leq\tau\leq\gamma$, $m\in\mathbb N_0$, we define the Sobolev norms
\begin{align*}
	\|\tens v\|_{\Omega_\tau}^m:= & \Bigl(\sum_{j=0}^m\int_{\Omega_\tau}|\widehat\nabla^j\tens v|^2\vect{\hat\mu}\Bigr)^{1/2}\text{ and} \\
	\|\tens v\|_{S_\tau}^m:= & \Bigl(\sum_{j=0}^m\int_{S_\tau}|\widehat\nabla^j\tens v|^2\vect{\hat\mu}_\tau\Bigr)^{1/2}.
\end{align*}
\end{Def}
Note that the $(n-1)$-dimensional Sobolev norm $\|\tens v\|_{S_\tau}^m$ is defined via the full $n$-dimensional derivative $\widehat\nabla$, i.\,e., derivatives are not restricted to the hypersurface $S_\tau$.
 	We emphasize that, in contrast to \cite{GMS09}, these Sobolev norms are completely of classical type, i.\,e., there is no $\ve$-dependence involved.
Equipped with this notion of Sobolev norms, we give the following definitions of \emph{energy tensors} and \emph{energy integrals}.
\begin{Def}
	For a smooth tensor field $\tens v$, $m>0$, and multiindices $K$ and $R$ with $|K|=|R|=m-1$, we define the \emph{energy tensors} $\tens T\ee^m(\tens v)$ of $\tens v$ of order $m$ by
	\begin{align*}
		T^{ab,0}\ee(\tens v):= & -\frac 1 2 g\ee^{ab}|\tens v|^2 \\
		T\ee^{ab,m}(\tens v):= & (g\ee^{ac}g\ee^{bd}-\frac 1 2 g\ee^{ab}g\ee^{cd}) m^{KR}m^{JQ}m_{IP}(\widehat\nabla_c\widehat\nabla_K v^I_J)(\widehat\nabla_d\widehat\nabla_{R} v^{P}_{Q}).
	\end{align*}
	For $0\leq\tau\leq\gamma$ and for $m\geq0$ we define the \emph{energy integral} $E^m_{\tau,\ve}(\tens v)$ of $\tens v$ of order $m$ on $S_\tau$ by
	\begin{equation*}
		E_{\tau,\ve}^m(\tens v):=\sum_{j=0}^m\int_{S_\tau}T\ee^{ab,j}(\tens v)\sigma_a\hat\sigma_b\vect{\hat\mu}_\tau.
	\end{equation*}
\end{Def}
To relate Sobolev norms and energy integrals, we need the following lemma, which is a variation of Lemma 4.1(1) in \cite{GMS09}.
\begin{Lemma}\label{L:Energies and Sobolev norms}
	There exist constants $A$ and $A'$ (independent of $\ve$) such that for all $k\geq 0$ and all smooth $\tens v$
	\begin{equation}\label{Eq:Energies and Sobolev norms}
		A'(\|\tens v\|^m_{S_\tau})^2 \leq E_{\tau,\ve}^m(\tens v)\leq A(\|\tens v\|^m_{S_\tau})^2.
	\end{equation}
	Consequently $E_{\tau,\ve}^m(\tens v)$ is $O(1)$.
\end{Lemma}
\begin{Proof}
	For $m=0$ we have
	\begin{align*}
		T\ee^{ab,0}(\tens v)\sigma_a\hat\sigma_b = & -\frac 1 2 g\ee^{ab}\sigma_a\hat\sigma_b|\tens v|^2=-\frac 1 2\cdot\frac 1 {\sigma} g\ee^{ab}\sigma_a\sigma_b|\tens v|^2 = -\frac{1}{2\sigma}\tens g\ee^{-1}(\vect\sigma,\vect\sigma)|\tens v|^2,
	\end{align*}
	which by \eqref{Eq:Bounds on g_e} yields
	\begin{equation*}
		A_0'|\tens v|^2\leq T\ee^{ab,0}(\tens v)\sigma_a\hat\sigma_b\leq A_0|\tens v|^2,
	\end{equation*}
	where we have set $A_0:=\frac{M_0}{2\sigma}$ and $A_0':=\frac 1{2\sigma M_0}$. So the result follows via integration over $S_\tau$.
	
	Now for $m>0$, we have
	\begin{align*}
		(g\ee^{ac}g\ee^{bd}-\frac 1 2g\ee^{ab}g\ee^{cd})\sigma_a\hat\sigma_b= & \frac 1{\sigma}(g\ee^{ac}g\ee^{bd}-\frac 1 2g\ee^{ab}g\ee^{cd})\sigma_a\sigma_b \\
		= & \frac 1{\sigma}(\xi^c\ee\xi\ee^d-\frac 1 2 \tens g\ee^{-1}(\vect\sigma,\vect\sigma) g\ee^{cd}) \\
		= & -\frac{\tens g\ee^{-1}(\vect\sigma,\vect\sigma)}{2\sigma}\left(-\frac{2\xi\ee^c\xi\ee^d}{\tens g\ee^{-1}(\vect\sigma,\vect\sigma)}+g\ee^{cd}\right).
	\end{align*}
	Hence, for $1\leq j\leq m$
	\begin{align*}
		T\ee^{ab,j}(\tens v)\sigma_a\hat\sigma_b= & -\frac{\tens g\ee^{-1}(\vect\sigma,\vect\sigma)}{2\sigma}\left(-\frac{2\xi\ee^c\xi\ee^d}{\tens g\ee^{-1}(\vect\sigma,\vect\sigma)}+g\ee^{cd}\right) \\
		& \cdot m^{KR}m^{JQ}m_{IP}(\widehat\nabla_c\widehat\nabla_Kv^I_J) (\widehat\nabla_d\widehat\nabla_{R}v^{P}_{Q}),
	\end{align*}
	where the expression in parentheses in the first line is a Riemannian metric since $\vect\xi\ee$ is timelike. Moreover, since we have that this Riemannian metric is locally bounded from above and below, we obtain the existence of $B$, $B'>0$ such that
	\begin{equation*}
		B'\tens m^{-1}(\vect\omega,\vect\omega)\leq \left(-\frac{2\vect\xi\ee\otimes\vect\xi\ee}{\tens g\ee^{-1}(\vect\sigma,\vect\sigma)}+\tens g\ee^{-1}\right)(\vect\omega,\vect\omega)\leq B\tens m^{-1}(\vect\omega,\vect\omega)
	\end{equation*}
	for any smooth one-form $\vect\omega$. Therefore
	\begin{equation*}
		-\frac{B'}{2\sigma}\tens g\ee^{-1}(\vect\sigma,\vect\sigma)|\widehat\nabla^j\tens u|^2\leq T\ee^{ab,j}(\tens u)\sigma_a\hat\sigma_b\leq -\frac{B}{2\sigma}\tens g\ee^{-1}(\vect\sigma,\vect\sigma)|\widehat\nabla^j\tens u|^2
	\end{equation*}
	and, again by \eqref{Eq:Bounds on g_e}, we have
	\begin{equation*}
		B'A'_0|\widehat\nabla^j\tens u|^2\leq T\ee^{ab,j}(\tens u)\sigma_a\hat\sigma_b\leq BA_0|\widehat\nabla^j\tens u|^2.
	\end{equation*}
	Finally, setting $A'=\min(A'_0,B'A'_0)$ and $A=\max(A_0,BA_0)$, integration, and summation over $j=0,\dots,m$ gives the result.
\end{Proof}

Next we note the essential fact that for all $\ve$, the energy tensors $\tens T\ee^j(\tens v)$ satisfy the \emph{dominant energy condition} with respect to $\tens g\ee$. That is, for any timelike one-form $\vect\omega\ee$, we have
$
	\tens T\ee^j(\tens v)(\vect\omega\ee,\vect\omega\ee)\geq 0,
$
and
$
	\tens T\ee^j(\tens v)(\vect\omega\ee,\cdot)$ is a non-spacelike vector field. This condition can also be entirely formulated in terms of $\mathcal G$, see \cite{May08}. 
	The dominant energy condition is a key feature in the following estimates: it guarantees positivity of
	\[
		\int_{\upartial\Omega_\tau\backslash(S_\tau\cup S_0)}\hspace{-1.8em}T\ee^{ab,j}(\tens v)\sigma_b\,\mathrm d \Omega_a,
	\]
	where $\mathrm d \Omega_a$ is the surface element on $\upartial\Omega$.
	Hence the dominant energy condition implies via the \emph{divergence theorem}, Lemma 4.3.1
	in \cite{HaEl73}, the following estimate:
\begin{align*}
	\int_{S_\tau} T\ee^{ab,j}(\tens v)\sigma_b\hat\sigma_a\vect{\hat\mu}_\tau\leq & \int_{S_0} T\ee^{ab,j}(\tens v)\sigma_b\hat\sigma_a\vect{\hat\mu}_0+C\int_0^\tau\!\!\int_{S_\zeta}T\ee^{ab,j}(\tens v)\sigma_b\hat\sigma_a\vect{\hat\mu}_\zeta\,\mathrm d\zeta \\
	& +\int_{\Omega_\tau}\widehat\nabla_a T\ee^{ab,j}(\tens v)\sigma_b\vect{\hat\mu}.
\end{align*}
Summation over $j$ for $0\leq\tau\leq\gamma$, yields the estimate
\begin{align}\label{Eq:Energy inequality 1}
	E_{\tau,\ve}^m(\tens v)\leq & E_{0,\ve}^m(\tens v)+\sum_{j=0}^m\biggl(C\int_0^\tau\!\!\int_{S_\zeta} T\ee^{ab,j}(\tens v)\sigma_a\hat\sigma_b\vect{\hat\mu}_\zeta\mathrm d \zeta + \int_{\Omega_\tau}\sigma_b\widehat\nabla_a T\ee^{ab,j}(\tens v)\vect{\hat\mu}\biggr)\notag \\
	= & E_{0,\ve}^m(\tens v) + C \int_0^\tau E_{\zeta,\ve}^m(\tens v)\mathrm d \zeta+\sum_{j=0}^m\,\int_{\Omega_\tau}\sigma_b\widehat\nabla_a T\ee^{ab,j}(\tens v)\vect{\hat\mu}
\end{align}
which will be our main tool in the following. In fact \eqref{Eq:Energy inequality 1} will be used to prove the energy estimates needed to derive moderateness and negligibility of solutions to \eqref{Eq:Representative PDE}.
The technical core is to provide estimates on the divergence term $\widehat\nabla_a T\ee^{ab,j}(\tens u\ee)$ in \eqref{Eq:Energy inequality 1} for a solution $\tens u\ee$ using the differential equation \eqref{Eq:Representative PDE}.
As mentioned in the introduction, in contrast to \cite{GMS09}, we deal with a connection with respect to the \emph{smooth} metric $\tens{\hat g}$, different from the ``coefficient metric'' $\tens g$. This amounts to additional terms containing $\widehat\nabla\tens g$, which we can control by condition (i) in Section 2.

\begin{Prop}\label{Prop:Energy inequality}
	Let $\tens u\ee$ be a solution of the differential equation \eqref{Eq:Representative PDE} on $U$. Then, for every $m\geq 1$, there exist constants $C_m'$, $C_m''$, and $C_m'''$ such that for every $0\leq\tau\leq\gamma$,
	\begin{align}\label{Eq:Energy inequality 2}
		E_{\tau,\ve}^m(\tens u\ee)\leq  & E_{0,\ve}^m(\tens u\ee)+C_m'(\|\tens F\ee\|_{\Omega_\tau}^{m-1})^2\notag \\
		& +C_m''\ve^{-N}
		\int_0^\tau E_{\zeta,\ve}^{m-1}(\tens u\ee)\,\mathrm d\zeta+C_m'''\int_0^\tau E_{\zeta,\ve}^m(\tens u\ee)\,\mathrm d\zeta.
	\end{align}
\end{Prop}
Observe that the coefficient in front of the last
integral in \eqref{Eq:Energy inequality 2} does not depend on $\ve$:
This is essential later on, when applying Gronwall's inequality in the
course of proving moderateness resp. negligibility of the nets $E_{\tau,\ve}^m(\tens u\ee)$.
\begin{Proof}
	To proof proposition \ref{Prop:Energy inequality}, we distinguish the cases $k=0$ and $k>0$. For $k=0$, we have
	\[
	\widehat\nabla_aT\ee^{ab,0}(\tens u\ee)=  -\frac 1 2\widehat\nabla_ag\ee^{ab}|\tens u\ee|^2-g\ee^{ab} u^I_{J,\ve}\widehat\nabla_a u^{P}_{Q,\ve}m^{JQ}m_{IP}
	 -\frac 1 2g\ee^{ab} u^I_{J,\ve} u^{P}_{Q,\ve}\widehat\nabla_a(m^{JQ}m_{IP}).
	\]
	Therefore, by the Cauchy-Schwartz inequality for the inner product induced by $\tens m$ on the tensor bundle on $M$, we have
	\begin{align*}
		\sigma_b\widehat\nabla_a T\ee^{ab,0}(\tens u\ee)\leq & \frac 1 2 |\sigma_b\widehat\nabla_a g\ee^{ab}|\cdot|\tens u\ee|^2+|\sigma_bg\ee^{ab}|\cdot|\widehat\nabla\tens u\ee|\cdot|\tens u\ee| \\
		& + \frac 1 2|\sigma_bg\ee^{ab}|\cdot|\widehat\nabla_a(m^{JQ}m_{IP})|\cdot|\tens u\ee|^2.
	\end{align*}
	This yields by condition (i) in Section \ref{Sec:Setting}
	\[
		\sigma_b\widehat\nabla_a T\ee^{ab,0}(\tens u\ee)\leq  P_{0}(|\tens u\ee|^2+|\widehat\nabla\tens u\ee|^2)
	\]
	for a constant $P_0$.
	
	For the case $k>0$, consider the expression
	\begin{align*}
		\sigma_b\widehat\nabla_aT\ee^{ab,k}(\tens u\ee)= & \sigma_b h\ee^{abcd}\widehat\nabla_a\widehat\nabla_c\widehat\nabla_K u^I_{J,\ve}\widehat\nabla_d\widehat\nabla_L u^{P}_{Q,\ve}m^{KL}m^{JQ}m_{IP} \\
		& +\sigma_b h\ee^{abcd}\widehat\nabla_c\widehat\nabla_K u^I_{J,\ve}\widehat\nabla_a\widehat\nabla_d\widehat\nabla_L u^{P}_{Q,\ve}m^{KL}m^{JQ}m_{IP} \\
		& +\sigma_b\widehat\nabla_a h\ee^{abcd}\widehat\nabla_c\widehat\nabla_K u^I_{J,\ve}\widehat\nabla_d\widehat\nabla_L u^{P}_{Q,\ve}m^{KL}m^{JQ}m_{IP} \\
		& +\sigma_b h\ee^{abcd}\widehat\nabla_c\widehat\nabla_K u^I_{J,\ve}\widehat\nabla_d\widehat\nabla_L u^{P}_{Q,\ve}\widehat\nabla_a (m^{KL}m^{JQ}m_{IP}),
	\end{align*}
	where $h\ee^{abcd}:=g\ee^{ac}g\ee^{bd}-\frac 12 g\ee^{ab}g\ee^{cd}$. After a rather lengthy calculation, where we interchange covariant derivatives, the above expression takes the form
	\begin{align}\label{Eq:Wahnsinnig komplizierter Ausdruck}
		\sigma_b\widehat\nabla_aT\ee^{ab,k}(\tens u\ee)= & \sigma_b g\ee^{bd}g\ee^{ac}\widehat\nabla_K\widehat\nabla_a\widehat\nabla_c u^I_{J,\ve}\widehat\nabla_d\widehat\nabla_L u^{P}_{Q,\ve}m^{KL}m^{JQ}m_{IP} \notag \\
		& + \sigma_b g\ee^{bd}g\ee^{ac}\Bigl(\sum_{j=0}^{k-1}(\mathcal R^{(k-1,j)}\tens u\ee)^I_{acKJ}\Bigr)\widehat\nabla_d\widehat\nabla_L u^{P}_{Q,\ve}m^{KL}m^{JQ}m_{IP} \notag \\
		& + \sigma_b g\ee^{bd}g\ee^{ac}(\mathcal R^{(k-1,k-1)}\tens u\ee)^I_{adKJ}\widehat\nabla_c\widehat\nabla_L u^{P}_{Q,\ve}m^{KL}m^{JQ}m_{IP} \notag \\
		& +\sigma_b\widehat\nabla_a h\ee^{abcd}\widehat\nabla_c\widehat\nabla_K u^I_{J,\ve}\widehat\nabla_d\widehat\nabla_L u^{P}_{Q,\ve}m^{KL}m^{JQ}m_{IP} \notag \\
		& +\sigma_b h\ee^{abcd}\widehat\nabla_c\widehat\nabla_K u^I_{J,\ve}\widehat\nabla_d\widehat\nabla_L u^{P}_{Q,\ve}\widehat\nabla_a (m^{KL}m^{JQ}m_{IP}).
	\end{align}
	Here the terms $\mathcal{\widehat R}^{(k,j)}\tens u\ee$ denote a linear combination of contractions of the $(k-j)$th derivative of $\tens{\hat R}$ with the $j$th derivative of $\tens u\ee$; $\tens {\hat R}$ being the Riemannian curvature tensor with respect to $\tens{\hat g}$.	
	In the next step, using the differential equation \eqref{Eq:Representative PDE} on the first line of \eqref{Eq:Wahnsinnig komplizierter Ausdruck}, we reduce the order of derivatives by one, obtaining $\sigma_b\widehat\nabla_aT\ee^{ab,k}(\tens u\ee)=\sum_{i=1}^{8}I_i$ with
	\begin{align*}
			I_1 = & \sigma_b g\ee^{bd}\widehat\nabla_K F^I_J\widehat\nabla_d\widehat\nabla_L u^{P}_{Q,\ve}m^{KL}m^{JQ}m_{IP}, \\
			I_2 = & -\sigma_b g\ee^{bd}\Bigl(\sum_{j=2}^k\mathcal A\ee^{(k+1,j)}\tens u\ee)^I_{KJ}\Bigr) \widehat\nabla_d\widehat\nabla_L u^{P}_{Q,\ve}m^{KL}m^{JQ}m_{IP}, \\
			I_3 = & -\sigma_b g\ee^{bd}\Bigl(\sum_{j=1}^k\mathcal B\ee^{(k,j)}\tens u\ee)^I_{KJ}\Bigr) \widehat\nabla_d\widehat\nabla_L u^{P}_{Q,\ve}m^{KL}m^{JQ}m_{IP}, \\
			I_4 = & -\sigma_b g\ee^{bd}\Bigl(\sum_{j=0}^{k-1}\mathcal C\ee^{(k-1,j)}\tens u\ee)^I_{KJ}\Bigr) \widehat\nabla_d\widehat\nabla_L u^{P}_{Q,\ve}m^{KL}m^{JQ}m_{IP}, \\
			I_5 = & \sigma_b g\ee^{bd}g\ee^{ac}\Bigl(\sum_{j=0}^{k-1}(\mathcal{\widehat R}^{(k-1,j)}\tens u\ee)^I_{acKJ}\Bigr)\widehat\nabla_d\widehat\nabla_L u^{P}_{Q,\ve}m^{KL}m^{JQ}m_{IP}, \\
			I_6 = & \sigma_b g\ee^{bd}g\ee^{ac}(\mathcal{\widehat R}^{(k-1,k-1)}\tens u\ee)^I_{adKJ}\widehat\nabla_c\widehat\nabla_L u^{P}_{Q,\ve}m^{KL}m^{JQ}m_{IP}, \\
			I_7 = & \sigma_b\widehat\nabla_a h\ee^{abcd}\widehat\nabla_c\widehat\nabla_K u^I_{J,\ve}\widehat\nabla_d\widehat\nabla_L u^{P}_{Q,\ve}m^{KL}m^{JQ}m_{IP}, \\
			I_8 = & \sigma_b h\ee^{abcd}\widehat\nabla_c\widehat\nabla_K u^I_{J,\ve}\widehat\nabla_d\widehat\nabla_L u^{P}_{Q,\ve}\widehat\nabla_a (m^{KL}m^{JQ}m_{IP}).
	\end{align*}
	Here $\mathcal A\ee^{(k,j)}\tens u\ee$, $\mathcal B\ee^{(k,j)}\tens u\ee$, and $\mathcal C\ee^{(k,j)}\tens u\ee$ denote a linear combination of contractions of the $(k-j)$th derivative of $\tens g\ee$, $\tens B\ee$, resp. $\tens C\ee$, resp. with the $j$th derivative of $\tens u\ee$.

	In the remaining part of the proof we study the $\ve$-asymptotics of the terms $I_1,\dots, I_8$. We aim at an estimate that involves no $\ve$-dependency in terms containing Sobolev norms of highest order. The following manipulations make repeatedly use of the Cauchy-Schwartz inequality. Proceeding in the same way as for the term of order $k=0$, we obtain
	\begin{align*}
		|I_1|\leq & P_{k,1} (|\widehat\nabla^k\tens u\ee|^2+|\widehat\nabla^{k-1}\tens F\ee|^2) \\
		|I_2|\leq & P_{k,2}\Bigl((k-1)|\widehat\nabla^k\tens u\ee|^2 + g_k^2\sum_{j=2}^{k-1}\ve^{-N_k}|\widehat\nabla^j\tens u\ee|^2\Bigr) \\
		|I_3|\leq & P_{k,3} \Bigl(k|\widehat\nabla^k\tens u\ee|^2+ b_k^2\sum_{j=1}^{k-1}\ve^{-N_k}|\widehat\nabla^j\tens u\ee|^2\Bigr) \\
		|I_4|\leq & \frac {P_{k,4}} 2 \Bigl(k|\widehat\nabla^k\tens u\ee|^2+ c_k^2\sum_{j=0}^{k-1}\ve^{-N_k}|\widehat\nabla^j\tens u\ee|^2\Bigr) \\
		|I_5|\leq & \frac {P_{k,5}} 2 \Bigl(k|\widehat\nabla^k\tens u\ee|^2+ r_k^2\sum_{j=0}^{k-1}|\widehat\nabla^j\tens u\ee|^2\Bigr) \\
		|I_6|\leq & \frac{P_{k,6}} 2 (|\widehat\nabla^k\tens u\ee|^2+r_k^2|\widehat\nabla^{k-1}\tens u\ee|^2) \\
		|I_7|\leq & P_{k,7}|\widehat\nabla^k\tens u\ee|^2
	\end{align*}
	\begin{align*}
		|I_8|\leq & P_{k,8}|\widehat\nabla^k\tens u\ee|^2
	\end{align*}
	for constants $P_{k,1},\dots,P_{k,8}$.
	Summing up, for $k>0$, there exist positive constants $\alpha_k$, $\beta_k$, and $\gamma_k$ such that
	\begin{equation*}
		\sigma_b\widehat\nabla_a T\ee^{ab,k}(\tens u\ee)\leq \alpha_k|\widehat\nabla^k\tens u\ee|^2+\beta_k|\widehat\nabla^{k-1}\tens F\ee|^2+\gamma_k\sum_{j=0}^{k-1}\ve^{-N_k}|\widehat\nabla^j\tens u\ee|^2.
	\end{equation*}
	Integration over $\Omega_\tau$ and summation over $k=0,\dots,m$ yields with $N:=\max_k N_k$
	\begin{align*}
		\sum_{k=0}^m\,\int_{\Omega_\tau}\sigma_b\widehat\nabla_a T\ee^{ab,k}(\tens u\ee)\hat\mu\leq & \tilde\alpha_m \int_0^\tau E^m_{\zeta,\ve}(\tens u\ee)\,\mathrm d\zeta+\tilde\beta_m(\|\tens F\ee\|_{\Omega_\tau}^{m-1})^2 \\
		& +\tilde\gamma_m\ve^{-N}
		\int_0^\tau E^{m-1}_{\zeta,\ve}(\tens u\ee)\,\mathrm d\zeta,
	\end{align*}
	where we used Lemma \ref{L:Energies and Sobolev norms} to translate Sobolev norms into energy integrals. Substitution of the last estimate into \eqref{Eq:Energy inequality 1} gives inequality \eqref{Eq:Energy inequality 2} and we are done.
\end{Proof}
Since we consider equations containing lower order terms, we need a sharper estimate for the first order energy integral. Thus we observe,
when setting $m=1$ in the proof of Proposition \ref{Prop:Energy inequality}, that the terms $I_2$ and $I_5$ vanish and the terms $I_1$, $I_6$, $I_7$, and $I_8$ do not contribute any $\ve$-powers. For the remaining terms $I_3$ and $I_4$, the coefficients $\tens B\ee$ and $\tens C\ee$ can be estimated by condition (i) in Section \ref{Sec:Setting}. Thus we have
\begin{Cor}\label{Cor:Sharper estimate}
	For $m=1$, we can sharpen inequality \eqref{Eq:Energy inequality 2} to
	\begin{equation}\label{Eq:Shaper energy inequality 2}
		E_{\tau,\ve}^1(\tens u\ee)\leq  E_{0,\ve}^1(\tens u\ee)+C_1'(\|\tens F\ee\|_{\Omega_\tau}^{0})^2 + C_1'''\int_0^\tau E_{\zeta,\ve}^1(\tens u\ee)\,\mathrm d\zeta.
	\end{equation}
\end{Cor}
\begin{Remark}\label{Rem:Compare with GMS09}
	Once more we undertake a comparison with \cite{GMS09}.
	\begin{enumerate}
		\item[1)] When looking at Proposition 5.1 in \cite{GMS09} and its proof one sees that the inclusion of lower order terms would give rise to an additional term for $j=0$ in the sum. This term, since dependent on $\ve$, would obstruct a successful application of Gronwall's inequality. To compensate for this $\ve$-dependence, we had to introduce bounds on $\tens B$ and $\tens C$ to obtain the vital Corollary \ref{Cor:Sharper estimate}.
		\item[2)] Since, in contrast to \cite{GMS09}, we do not work with the Levi-Cività connection of the ``coefficient metric'', we had to sharpen the condition on first order derivatives of $\tens g\ee$ to $\sup_K|\widehat\nabla \tens g^{-1}\ee|=O(1)$, see (i).
		\item[3)] As the present analysis shows, the assumptions in \cite{GMS09}, condition (A) on the derivatives of $\tens g$ are not necessary to prove the existence and uniqueness result \cite[Thm. 3.1]{GMS09}. In fact, one may significantly improve the main theorem in \cite{GMS09} by relaxing (A) to $\sup_K |\tens g\ee|=O(1)$ and $\sup_K |\tens g^{-1}\ee|=O(1)$. On the other hand, the original condition (A) obviously can be used to derive the precise asymptotics of the solutions of \cite[Thm. 3.1]{GMS09}, hence to prove an additional regularity result.
	\end{enumerate}
\end{Remark}

We may now apply Gronwall's inequality to \eqref{Eq:Energy inequality 2} and \eqref{Eq:Shaper energy inequality 2} to immediately obtain
\begin{Cor}\label{Cor:Energy inequality}
	Let $\tens u\ee$ be a solution of the differential equation \eqref{Eq:Representative PDE} on $U$. Then, for every $m> 1$, there exist constants $C_m'$, $C_m''$, and $C_m'''$ such that for every $0\leq\tau\leq\gamma$,
	\begin{align}\label{Eq:Energy inequality 3}
		E_{\tau,\ve}^m(\tens u\ee)\leq  & \biggl(E_{0,\ve}^m(\tens u\ee)+C_m'(\|\tens F\ee\|_{\Omega_\tau}^{m-1})^2+C_m''\ve^{-N}
		\int_0^\tau E_{\zeta,\ve}^{m-1}(\tens u\ee)\,\mathrm d\zeta\biggr)e^{C_m'''\tau}.
	\end{align}
	For $m=1$, we have
	\begin{align}\label{Eq:Energy inequality 4}
		E_{\tau,\ve}^1(\tens u\ee)\leq  & \biggl(E_{0,\ve}^1(\tens u\ee)+C_1'(\|\tens F\ee\|_{\Omega_\tau}^{0})^2\biggr)e^{C_1'''\tau}.
	\end{align}
\end{Cor}
The consequence of \eqref{Eq:Energy inequality 3} and \eqref{Eq:Energy inequality 4} is that by iterating $m$, we obtain that moderate resp. negligible initial energy integrals $E^m_{0,\ve}(\tens u\ee)$ and right hand side $\tens F\ee$ imply moderate resp. negligible energy integrals $E^m_{\tau,\ve}(\tens u\ee)$ at later times $\tau$. We will make use of this fact in the proof of the main theorem.
\begin{Cor}\label{Cor:Moderate Energy}
	Let $\tens u\ee$ be a solution of the initial value problem \eqref{Eq:Representative PDE} on $U$. Then for all $m\geq 1$, if the initial energy integrals are moderate,
	\begin{equation}\label{Eq:Moderate Energy}
		\sup_{0\leq\tau\leq\gamma}E_{\tau,\ve}^m(\tens u\ee)\ee
	\end{equation}
	 is a moderate net of real numbers. Likewise, \eqref{Eq:Moderate Energy} is negligible if the initial energy integrals and, additionally, $\tens F$ is negligible.
\end{Cor}

\section{Proof of the main theorem}

Before we start the actual proof of the main theorem we need two more estimates. First, we translate bounds on initial data into bounds on initial energy integrals and bounds on energy integrals into bounds on solutions. Afterwards, together with Corollary \ref{Cor:Moderate Energy}, we will establish existence and uniqueness of generalized solutions to the initial value problem \eqref{Eq:Wave equation}, \eqref{Eq:Initial data}.

\begin{Lemma}\label{L:Initial energies from initial data}
	Let $(\tens u\ee)\ee$ be a solution of \eqref{Eq:Representative PDE}. If $(\tens u_{0,\ve})\ee$, $(\tens u_{1,\ve})\ee$ and $(\tens F\ee)\ee$ are moderate resp. negligible, then the initial energy integrals $(E^m_{0,\ve}(\tens u\ee))\ee$ for each $m\geq 0$ are moderate resp. negligible nets of real numbers.
\end{Lemma}
\begin{Proof}
	We have to consider the energy integral
	\begin{equation*}
		E_{0,\ve}^m(\tens u\ee)=\sum_{j=0}^m\int_{S_0} T\ee^{ab}(\tens u\ee)\sigma_a\hat\sigma_b\vect\mu_\tau.
	\end{equation*}
	Obviously one can see that the moderateness resp. negligibility of the initial energy integrals $E_{0,\ve}^1(\tens u\ee)$ is equivalent to the moderateness resp. negligibility of $\tens u\ee$ and its first order derivatives on $S_0$, thus immediately follows from moderatenes resp. negligibility of the data. To deal with the higher order initial energy integrals, we need further arguments: 
	Choosing a coordinate system $(t,x^\alpha)=(x^0,x^\alpha)$, we have
	\begin{align*}
		\tens u\ee(0,x^\alpha)= & \tens u_{0,\ve}(x^\alpha) \\
		\upartial_t\tens u\ee(0,x^\alpha)= & \tens{\tilde u}_{1,\ve}(x^\alpha),
	\end{align*}
	where
	\begin{equation}\label{Eq:Initial condition in coordinates}
		\tens{\tilde u}_{1,\ve}:=\frac 1{\xi\ee^0}\Bigl(\tens u_{1,\ve}-\sum_{\lambda=1}^n\xi\ee^\lambda\widehat\nabla_\lambda\tens u_{0,\ve}-\xi\ee^0(\widehat\nabla_t-\upartial_t)\tens u_{0,\ve}\Bigr)
	\end{equation}
	 and $(\widehat\nabla_t-\upartial_t)\tens u_{0,\ve}$ denotes the difference between the covariant and the partial derivative of $\tens u_{0,\ve}$, which may be expressed in terms of the Christoffel symbols $\tens{\hat\Gamma}$.
Now terms of the energy tensors including first order time-derivatives are by \eqref{Eq:Initial condition in coordinates} rewritten in terms of the (known) spatial derivatives of the data. For higher order time-derivatives, we inductively use the differential equation in the form
	 \begin{align*}
	 	\upartial^2_t u^{\mu_1\cdots\mu_k}_{\nu_1\cdots\nu_l,\ve}= & -\frac 1{g\ee^{00}}\Bigl(\sum_{\mu,\nu=1}^n g\ee^{\mu\nu}\upartial_\mu\upartial_\nu u^{\mu_1\cdots\mu_k}_{\nu_1\cdots\nu_l,\ve} \\
	& +\text{terms with less than 2 time derivatives}\Bigr)
	 \end{align*}
	to reduce the order of time-derivatives on $\tens u$ to 1.
	By assumption all the coefficients occurring in the differential equation are moderate resp. negligible, thus not contributing more than a factor $\ve^{-N}$ for some natural number $N$.
\end{Proof}
\begin{Lemma}\label{L:Bounds on solutions from energies}
	Let $\alpha$ be a multiindex with $|\alpha|=m$. For $s>(n-1)/2$ an integer, there exists a constant $C$ and number $N$ such that for all $\tens u\in\mathcal T^k_l(\Omega_\tau)$ and for all $\zeta\in[0,\tau]$, we have
	\begin{equation*}
		\sup_{p\in\Omega_\tau}|\upartial^\alpha\tens u(x)|\leq \frac C{\sqrt{A'}}
		\sup_{0\leq\zeta\leq\tau} (E^{s+m}_{\zeta,\ve}(\tens u))^{1/2}.
	\end{equation*}
\end{Lemma}
\begin{Proof}
	By the Sobolev embedding theorem on $S_\tau$
	we obtain for $s>(n-1)/2$,
	\begin{equation}\label{Eq:Sobolevenergies}
		\sup_{p\in S_\zeta}|\tens u(x)|\leq C  \|\tens u\|^s_{S_\zeta}.
	\end{equation}
	Application of \eqref{Eq:Energies and Sobolev norms} yields
	\begin{equation*}
		\sup_{p\in S_\zeta}|\tens u(x)|\leq \frac{C}{\sqrt {A'}} (E^s_{\zeta,\ve}(\tens u))^{1/2}.
	\end{equation*}
	We then take the supremum over $\zeta\in[0,\tau]$ on the right hand side and obtain the result for $m=0$. To show the general result, we replace $\tens u$ by the respective derivatives, i.\,e., we replace \eqref{Eq:Sobolevenergies} by
	\begin{equation*}
		\sup_{p\in S_\zeta}|\upartial_{\rho_1}\cdots\upartial_{\rho_i}\upartial_t^{j} \tens u|\leq C \|\upartial_t^{j}\tens u\|_{S_\zeta}^{s+i}\leq C\|\tens u\|_{S_\zeta}^{s+i+j},
	\end{equation*}
	where $m=i+j$
\end{Proof}
Finally, we are ready to proof Theorem \ref{Th:Main Theorem}.

\vspace{0.5\baselineskip}
\noindent\textbf{Proof of the main theorem:}
	
\vspace{0.5\baselineskip}
\noindent
\emph{Step 1: Existence of classical solutions}. \\ Theorem 5.3.2 in \cite{Fri75} together with assumption (ii) in Section \ref{Sec:Setting} guarantees existence of a unique smooth solution $\tens u\ee$ of \eqref{Eq:Representative PDE} for each $\ve$ on a domain $A\subset\bigcap_{\ve<\ve_0} \overline{I\ee^+(\Sigma)}$. Without loss of generality we may assume, that $\Omega_\gamma\subset A$.

\vspace{0.5\baselineskip}
\noindent\emph{Step 2: Existence of generalized solutions}. \\ We show that the net obtained in step 1 is moderate on $\Omega_\gamma$. By assumption $\tens u_0$ and $\tens u_1$ are moderate, so by Lemma \ref{L:Initial energies from initial data} we obtain moderate initial energy integrals $E_{0,\ve}^k(\tens u\ee)$ for $k\geq 1$ and $0\leq\tau\leq\gamma$. Now Corollary \ref{Cor:Moderate Energy} ensures moderateness of the energy integrals $E_{\tau,\ve}^k(\tens u\ee)$ for $k\geq 1$ and $0\leq\tau\leq\gamma$. Finally, Lemma \ref{L:Bounds on solutions from energies} implies moderateness of $\tens u\ee$, hence $\tens u:=[(\tens u\ee)\ee]$ is a generalized solution of the initial value problem \eqref{Eq:Wave equation}, \eqref{Eq:Initial data} on $\Omega_\gamma$.

\vspace{0.5\baselineskip}
\noindent\emph{Step 3: Independence of the representatives of the data}. \\ The proof follows the same arguments as in step 2. Since $L$ is a linear differential operator, it suffices to show that the solution $\tens u$ of an equation with negligible $\tens F$, $\tens u_0$, and $\tens u_1$ is negligible as well. To establish this result, we proceed as before using the negligibility parts of Lemma \ref{L:Initial energies from initial data} and Corollary \ref{Cor:Moderate Energy}.
Thus the solution $\tens u=[(\tens u\ee)\ee]$ is unique and we are done.\hfill$\Box$

\section*{Acknowledgements}

The author wants to express his gratitude to the organizers of the \emph{International Conference of Generalized Functions 2009}, especially the Viennese part of the DIANA research group. This work was supported by research grants P20525 and Y237 of the Austrian Science Fund and the Initiativkolleg 1008-N \emph{Differential Geometry and Lie Groups} of the University of Vienna.

\end{document}